\documentclass{elsart}%
\usepackage{amsmath}
\usepackage{graphicx}
\usepackage{enumerate}
\usepackage{float}
\usepackage{indentfirst}
\usepackage{verbatim}
\usepackage{algorithm}
\usepackage{algpseudocode}
\usepackage{float}
\usepackage{multicol}
\usepackage{listings}
\usepackage{subfigure}
\usepackage{epsfig}
\usepackage{amsfonts}
\usepackage{amssymb}
\usepackage{epstopdf}%
\everymath{\displaystyle}
\begin{document}
\begin{frontmatter}
\title{Numerical Solution of the Robin Problem of Laplace Equations with
a Feynman-Kac Formula and Reflecting Brownian Motions}
\author[UNCC]{Yijing Zhou},
\author[UNCC]{Wei Cai}
\address[UNCC]{Department of Mathematics and Statistics,
University of North Carolina at Charlotte, Charlotte, NC 28223-0001}
\bigskip
\newpage
\begin{abstract}
In this paper, we present numerical methods to implement the probabilistic representation
of third kind (Robin) boundary problem for the Laplace equations.
The solution is based on a Feynman-Kac formula for the Robin problem which
employs the standard reflecting Brownian motion
(SRBM) and its boundary local time arising from the Skorohod problem.
By simulating SRBM paths through Brownian motion using Walk on Spheres
(WOS) method, approximation of the boundary local time is obtained and the Feynman-Kac formula
is calculated by evaluating the average of all path integrals over the boundary under a measure
defined through the local time. Numerical results demonstrate the accuracy and efficiency of the
proposed method for finding a local solution of the Laplace equations with Robin boundary conditions.
\end{abstract}
\end{frontmatter}
\numberwithin{equation}{section}

\section{Introduction}

Partial differential equations (PDEs) have been widely used to describe a
variety of phenomena such as electrostatics, electrodynamics, fluid flow or
quantum mechanics. Traditionally, finite difference, finite element and
boundary element methods are the mainstream numerical approaches to solve the
PDEs. Recently, using the Feynman-Kac formula \cite{[5]}\cite{[6]}\cite{[7]}
which connects solutions of differential equations of diffusion and heat flow
and random processes of Brownian motions, numerical methods based on random
walks or Monte Carlo diffusions have been explored for solving parabolic and
elliptic PDEs \cite{[11]}\cite{[24]}.

The Feynman-Kac formula represents the solutions of parabolic and elliptic
PDEs as the expectation functionals of stochastic processes (specifically
Brownian motions), and conversely, the probabilistic properties of diffusion
processes can be obtained through investigating related PDEs characterized by
corresponding generators \cite{[21]}. The formula involves the path integrals
of the diffusion process starting from an arbitrarily prescribed location, and
this enables us to find a local numerical solution without constructing space
and time meshes as in traditional deterministic numerical methods mentioned
above, which incur expensive costs in high dimensions. In many applications it
is also of practical importance and necessity to seek a local solution of PDEs
at some interested points. If the sample paths of a diffusion process are
simulated, then by computing the average of path integrals we can obtain
approximations to the exact solutions of the PDEs. For second order elliptic
PDEs with Dirichlet and Neumann boundaries, the average of path integrals is
reduced to the average of boundary integrals under certain measure where the
detailed trajectories of the diffusion process have no effect on the averages
except the hitting locations on the boundaries.

Simulations of diffusion paths can be done by random walks methods
\cite{[3]}\cite{[8]} \cite{[11]} \cite{[13]} either on lattice or in continuum
space. In some cases such as for the Poisson equation, the Feynman-Kac formula
has a pathwise integral requiring the detailed trajectory of each path.
Moreover, one may need to adopt random walks on a discrete lattice in order to
incorporate inhomogeneous source terms. As for the continuum space approach,
the Walk on Spheres (WOS) method is preferred where the path of diffusion
process within the domain does not appear in the Feynman-Kac formula. For both
approaches, the geometry of the boundaries need special care for accurate
results \cite{[14]}. In our previous work on Laplace equation with Neumann
boundary conditions \cite{[4]}, we proposed a numerical method to simulate the
standard reflecting Brownian motion (SRBM) path using WOS and obtained the
boundary local time of the SRBM. As a result, a local numerical solution of
the PDE is achieved by using the Feynman-Kac formula. Other literatures
\cite{[9]}\cite{[10]}\cite{[13]}\cite{[14]} have also explored similar
problems. Especially, in \cite{[13]} schemes based on the WOS, Euler schemes
and kinetic approximations are proposed to treat inhomogeneous Neumann
problems. It turns out that the pointwise resolution is much harder due to the
choice of the truncation of time. However, the local time was not handled
explicitly in \cite{[13]}. On the other hand, Monte Carlo simulations were
discussed in \cite{[14]} where the positive part of the boundary needs to be
identified first. In this paper, following \cite{[4]} we continue the use of
SRBM to solve Robin boundary problems for the Laplace operator, which has many
applications in heat transfer and impedance tomograph. Our goal again is to
obtain a local approximation to the exact solution of the Robin problem.

The rest of paper is organized as follows. Firstly, the Skorohod problem is
introduced in section 2, where both the concepts of standard reflecting
Brownian motion and boundary local time will be reviewed briefly. This lays
the foundation for the underlying diffusion process of the Robin boundary
problem and the sampling of the diffusion paths. Secondly, an overview of the
Feynamn-Kac formula is given in section 3. Thirdly, the probabilistic
representation of the solution for the Robin boundary value problem proposed
in \cite{[2]}\cite{[12]} is discussed in section 4, and we will see the
relation between the Neumann and Robin problems and gain a new perspective.
Section 5 presents the numerical approaches and test results. Finally,
conclusions and future work are given in section 6.

\section{Skorohod problem, SRBM and boundary local time}

Assume that $D$ is a domain with a $C^{1}$ boundary in $R^{3}$. The
generalized Skorohod problem is stated as follows: \newline

\begin{defn}
Let $f\in C([0,\infty),R^{3})$, a continuous function from $[0,\infty]$ to
$R^{3}$. A pair $(\xi_{t},L_{t})$ is a solution to the Skorohod equation
$S(f;D)$ if

\begin{enumerate}
\item $\xi$ is continuous in $\bar{D}$;

\item $L(t)$ is a nondecreasing function which increases only when $\xi
\in\partial D$, namely,
\begin{equation}
L(t)=\int_{0}^{t}I_{\partial D}(\xi(s))L(ds); \label{eq1}%
\end{equation}

\item The Skorohod equation holds:
\begin{equation}
S(f;D):\qquad\ \xi(t)=f(t)-\frac{1}{2}\int_{0}^{t}n(\xi(s))L(ds), \label{eq3}%
\end{equation}
where $n(x)$ denotes the outward unit normal vector at $x\in\partial D$.
\end{enumerate}
\end{defn}

The Skorohod problem was first studied in \cite{[1]} by A.V. Skorohod in
addressing the construction of paths for diffusion processes with boundaries,
which results from the instantaneous reflection behavior of the processes at
the boundaries. Skorohod presented the result in one dimension in the form of
an Ito integral and Hsu \cite{[12]} later extended the concept to
$d$-dimensions ($d\geq2$).

In the simple case that $D=[0,\infty)$, the solution to the Skorohod problem
uniquely exists and can be explicitly given by
\begin{equation}
\centering\xi(t)=\left\{
\begin{aligned} &f(t), &if\  t\leq\tau;\\ &f(t)-\inf_{\tau\leq s\leq t}f(s), &if\  t>\tau;\\ \end{aligned}\right.
\end{equation}
where $\tau=\inf\left\{  t>0:f(t)<0\right\}  $. In general, solvability of the
Skorohod problem is closely related to the smoothness of the domain $D$. For
higher dimensions, the existence of ($\ref{eq3}$) is guranteed for $C^{1}$
domains while uniqueness can be acheived for a $C^{2}$ domain by assuming the
convexity for the domain \cite{[15]}. Later, it was shown by Lions and
Sznitman \cite{[16]} that the constraints on $D$ can be relaxed to some
locally convex properties.

Next we introduce the concept of SRBM and boundary local time which play
important roles in solving Robin boundary problem by probabilistic approaches.

Suppose that $f(t)$ is a standard Brownian motion (SBM) starting at $x\in
\bar{D}$ and $(X_{t},L_{t})$ is the solution to the Skorohod problem $S(f;D)$,
then $X_{t}$ will be the standard reflecting Brownian motion (SRBM) on $D$
starting at $x$. Because the transition probability density of the SRBM
satisfies the same parabolic differential equation as that by a BM, a sample
path of the SRBM can be simulated simply as that of \ the BM within the
domain. However, the zero Neumann boundary condition for the density of SRBM
implies that the path be pushed back at the boundary along the inward normal
direction whenever it attempts to cross the latter. The full construction of a
SRBM from a SBM can be found in our previous work \cite{[4]}.

The boundary local time $L_{t}$ is not an independent process but associated
with SRBM $X_{t}$ and defined by
\begin{equation}
L(t)\equiv\lim_{\epsilon\rightarrow0}\frac{\int_{0}^{t}I_{D_{\epsilon}}%
(X_{s})ds}{\epsilon}, \label{eq5}%
\end{equation}
where $D_{\epsilon}$ is a strip region of width $\epsilon$ containing
$\partial D$ and $D_{\epsilon}\subset\overline{D}$. Here $L_{t}$ is called the
local time of $X_{t}$, a notion invented by P. L\'{e}vy \cite{[22]}. This
limit exists both in $L^{2}$ and $P^{x}$-$a.s$. for any $x\in\overline{D}$.

It is obvious that $L_{t}$ measures the amount of time that the standard
reflecting Brownian motion $X_{t}$ spends in a vanishing neighborhood of the
boundary within the time period $[0,t]$. Besides, it is the unique continuous
nondecreasing process that appears in the Skorohod equation. An interesting
part of ($\ref{eq5}$) is that the set $\left\{  t\in R_{+}:X_{t}\in\partial
D\right\}  $ has a zero Lebesgue measure while the sojourn time of the set is
nontrivial \cite{[22]}. This concept is not just a mathematical one but also
has physical relevance in understanding the \textquotedblleft crossover
exponent" associated with \textquotedblleft renewal rate" in modern renewal
theory \cite{[17]}.

In \cite{[12]}, an alternative explicit form of the local time was found,
\begin{equation}
L(t)=\sqrt{\frac{\pi}{2}}\int_{0}^{t}I_{\partial D}(X_{s})\sqrt{ds},
\label{eq7}%
\end{equation}
where the the right-hand side of (\ref{eq7}) is understood as the limit of
\begin{equation}
\sum_{i=1}^{n-1}%
\smash{\displaystyle\max_{s\in\Delta_i}I_{\partial D}(X_s)\sqrt{|\Delta_i|}},\quad
\smash{\displaystyle\max_i}|\Delta_{i}|\rightarrow0, \label{eq9}%
\end{equation}
where $\Delta=\{\Delta_{i}\}$ is a partition of the interval $[0,t]$ and each
$\Delta_{i}$ is an element in $\Delta$. ($\ref{eq5}$) and ($\ref{eq7}$)
provide us different ways to approximate local time and in \cite{[4]}, it was
found that ($\ref{eq5}$) yields better approximations in Neumann problem than
($\ref{eq7}$). Therefore, in this paper, we will also choose ($\ref{eq5}$) as
the approach to estimate the local time here.

\section{A Feynman-Kac formula}

The Feynman-Kac formula named after Richard Feynman and Mark Kac, establishes
a link between PDEs and stochastic processes. It first arose in the potential
theory for Sch\"{o}dinger equations, leading to a profound reformulation of
the quantum mechanics by the means of path integrals. Later, the formula also
finds its applications in mathematical finance, where the probabilistic and
the PDE representations in derivative pricing are connected.

Let us first look at the Dirichlet problems. Given a domain $D\subset R^{d}$
with a boundary $\partial D$,
\begin{equation}
\left\{
\begin{aligned} Lu(x)-c(x)u(x)&=f(x), \ x\in D\\ u(x)&=\phi (x), \ x\in\partial D\\ \end{aligned}\right.
, \label{eq11}%
\end{equation}
where the operator $L=-\frac{1}{2}\sum_{i,j=1}^{d}a_{ij}(x)\frac{\partial^{2}%
}{\partial x^{i}\partial x^{j}}-\sum_{i=1}^{d}b_{i}(x)\frac{\partial}{\partial
x^{i}}$ and both the coefficients in $L$ and $c(x)$ are Lipschitz continuous
and bounded.

The Feynman-Kac formula in this case \cite{[18]} represents the solution to
($\ref{eq11}$) in terms of an Ito diffusion process $X_{t}(\omega),$
\begin{equation}
u(x)=E^{x}[\int_{0}^{\tau_{D}}f(X_{t})exp\left\{  \int_{0}^{t}c(X_{s}%
)ds\right\}  dt]+E^{x}[\phi(X_{\tau_{D}})exp\left\{  \int_{0}^{\tau_{D}%
}c(X_{s})ds\right\}  ], \label{eq13}%
\end{equation}
with $\tau_{D}=\inf\{t:X_{t}\in\partial D\}$ and $X_{t}(\omega)$ is defined
by
\begin{equation}
dX_{t}=b(X_{t})dt+\alpha(X_{t})dB_{t}, \label{eq15}%
\end{equation}
where $B_{t}$ is the Brownian motion and $[a_{ij}]=\frac{1}{2}\alpha
(x)\alpha^{T}(x),[b_{ij}]=b$.

The expectation $E^{x}$ is an integration with respect to a measure $P_{x}$
taken over all sample paths $X_{t=0}(\omega)=x$, thus ($\ref{eq13}$) is a
representation of a solution of Dirichlet problem in the form of functional
integral. Moreover, ($\ref{eq13}$) is obtained by killing process $X_{t}$ at a
stopping time $\tau_{D}$ at which $X_{t}$ will be absorbed on the boundary. If
$c(x)\geq0$, then the function $c(x)$ can be interpreted as the killing rate
\cite{[21]}. It should be pointed out that (\ref{eq13}) is equivalent to the
formulation of weak solution and it is a classical solution as well if some
smoothness conditions are satisfied.

The Feynman-Kac formula above offers a method for solving certain PDEs by
simulating random paths of a stochastic process. Conversely, an important
class of expectations of random processes can be computed by deterministic
methods. For the Neumann boundary condition, a similar formula was derived by
Hsu \cite{[12]} for the Poisson equation, which is in the form of a functional
integral based on the boundary local time introduced in section 2. In this
case, though the Feynman-Kac formula remains in a similar form, it should be
understood as a path integral over the stochastic process $L_{t}$ associated
with the standard reflecting Brownian motion.

\section{Robin boundary value problem}

We focus on Robin boundary value problem for the time-independent
Schr\"{o}dinger equation.
\begin{equation}
\centering\left\{
\begin{aligned} \frac{1}{2}\Delta u+qu&=0,\quad in\  D;\\ \frac{\partial u}{\partial n}-cu&=f,\quad on\  \partial D.\\ \end{aligned}\right.
\label{eq17}%
\end{equation}

A generalization of the Feynman-Kac formula of section 3 in \cite{[2]} gives a
probablistic solution of ($\ref{eq17}$) as follows,
\begin{equation}
u(x)=E^{x}\left\{  \int_{0}^{\infty}e_{q}(t)\hat{e}_{c}(t)f(X_{t}%
)dL_{t}\right\}  , \label{eq19}%
\end{equation}
where $X_{t}$ is a SRBM starting at $x$. The term Feynman-Kac functional
$e_{q}(t),$ also appeared in the Neumann problem \cite{[12]}, is defined as
\begin{equation}
e_{q}(t)=\exp\left[  \int_{0}^{t}q(X_{s})\,ds\right]  , \label{eq21}%
\end{equation}
and a second functional is introduced for the Robin boundary problem, for
$c\in\Sigma_{d}(\partial D)$
\begin{equation}
\hat{e}_{c}(t)=exp\left[  \int_{0}^{t}c(X_{s})dL_{s}\right]  . \label{eq23}%
\end{equation}

Using these two functionals, we have,
\begin{equation}
u(x)=E^{x}\left\{  \int_{0}^{\infty}exp\left[  \int_{0}^{t}\left(
q(X_{s})ds+c(X_{s})dL_{s}\right)  \right]  f(X_{t})dL_{t}\right\}  .
\label{eq25}%
\end{equation}

Recalling the definition of the local time in ($\ref{eq5}$), we have the
following approximation
\begin{equation}
L(t)\approx\frac{1}{\epsilon}\int_{0}^{t}I_{D_{\epsilon}}(X_{s})ds,
\label{eq27}%
\end{equation}
thus,
\begin{equation}
dL(s)\approx\frac{1}{\epsilon}I_{D_{\epsilon}}(X_{s})ds. \label{eq29}%
\end{equation}

Therefore, ($\ref{eq25}$) can be modified as
\begin{equation}
u(x)\approx E^{x}\left\{  \int_{0}^{\infty}exp\left[  \int_{0}^{t}\left(
q(X_{s})+\frac{1}{\epsilon}c(X_{s})I_{D_{\epsilon}}(X_{s})\right)  ds\right]
f(X_{t})dL_{t}\right\}  , \label{eq31}%
\end{equation}
It can also be shown that as $\epsilon$ goes to zero, ($\ref{eq31}$) converges
to ($\ref{eq25}$) uniformly on $\bar{D}$.

As (\ref{eq31}) resembles the Feyman-Kac formula for the Neumann problem with
a modified $q(x)$ \cite{[4]}, it indicates a connection between the Robin and
the Neumann problems, namely, we may introduce
\begin{equation}
q_{\epsilon}(x)=q(x)+\frac{1}{\epsilon}c(x)I_{D_{\epsilon}}(x), \label{eq33}%
\end{equation}
then, the Robin boundary problem ($\ref{eq17}$) can be viewed as a limiting
case ($\epsilon\rightarrow0$) of Neumann problems
\begin{equation}
\centering\left\{
\begin{aligned} \frac{1}{2}\Delta u+q_\epsilon u&=0,\quad in\  D;\\ \frac{\partial u}{\partial n}&=f,\quad on\  \partial D.\\ \end{aligned}\right.
\label{eq35}%
\end{equation}

\section{Numerical approach and results}

In the present work, we only consider the case of the Laplace equation where
$q=0$ in ($\ref{eq35}$). From ($\ref{eq25}$),
\begin{equation}
u(x)=E^{x}\left\{  \int_{0}^{\infty}e^{\int_{0}^{t}c(X_{t})dL_{t}}%
f(X_{t})dL_{t}\right\}  , \label{eq37}%
\end{equation}
where $X_{t}$ represents the standard reflecting Brownian motion. For the sake
of computer simulation, the time period is truncated into $[0,T]$ to produce
an approximation for $u(x)$, i.e.,
\begin{equation}
\tilde{u}(x)=E^{x}\left\{  \int_{0}^{T}e^{\int_{0}^{t}c(X_{t})dL_{t}}%
f(X_{t})dL_{t}\right\}  . \label{eq39}%
\end{equation}

Next we will give a general description on the realization of SRBM paths and
the calculation of the corresponding local time, as implemented in \cite{[4]}.
A SRBM path can be constructed by pulling back a BM path back onto the
boundary whenever it runs out of the domain. Specifically, a SRBM path behaves
exactly the same way as a BM which is simulated by the WOS method.

\subsection{Simulating SRBM by the method of Walk on Spheres (WOS)}

\begin{itemize}
\item Method of WOS for Brownian paths
\end{itemize}

Random walk on spheres (WOS) method was first proposed by M\"{u}ller
\cite{[7]}, which can solve the Dirichlet problem for the Laplace operator
efficiently \cite{[8]}\cite{[10]} .

To illustrate the WOS method for the Dirichlet problem (\ref{eq11}), let us
consider the Laplace equation again where $f=0,a_{ij}=\delta_{ij}$ and
$b_{i}=0$ in (\ref{eq11}) and the It\^{o} diffusion is then simply the
standard Brownian motion with no drift. The solution to the Laplace equation
can be rewritten in terms of a measure $\mu_{D}^{x}$ defined on the boundary
$\partial D$,
\begin{equation}
u(x)=E^{x}(\phi(X_{\tau_{D}}))=\int_{\partial D}\phi(y)d\mu_{D}^{x},
\end{equation}
where $\mu_{D}^{x}$ is the harmonic measure\ defined by
\begin{equation}
\mu_{D}^{x}(F)=P^{x}\left\{  X_{\tau_{D}}\in F\right\}  ,F\subset\partial
D,x\in D.
\end{equation}
It can be shown easily that the harmonic measure is related to the Green's
function $g(y,x)$ for the domain with a homogeneous boundary condition
\cite{CKL}, i.e.,
\begin{equation}
\left\{
\begin{aligned} -\Delta g(x,y) &= \delta(x-y), \  &x\in D,\\ g(x,y) &= 0, \  &x\in\partial D\\ \end{aligned}\right.
,
\end{equation}
as follows
\begin{equation}
p(\mathbf{x},\mathbf{y})=-\frac{\partial g(x,y)}{\partial n_{y}}.
\end{equation}

If the starting point $x$ of a Brownian motion is at the center of a ball, the
probability of the BM exiting a portion of the boundary of the ball will be
proportional to the portion's area. Therefore, sampling a Brownian path by
drawing balls within the domain can significantly reduce the path sampling
time. To be specific, given a starting point $x$ inside the domain $D$, we
simply draw a ball of largest possible radius fully contained in $D$ and then
the next location of the Brownian path on the surface of the ball can be
sampled, using a uniform distribution on the sphere, say at $x_{1}$. Treat
$x_{1}$ as the new starting point, draw a second ball fully contained in $D$,
make a jump from $x_{1}$ to $x_{2}$ on the surface of the second ball as
before. Repeat this procedure until the path hits a absorption $\epsilon
$-shell of the domain (see Fig. 2) \cite{[5]}. When this happens, we assume
that the path has hit the boundary $\partial D$ (see Fig. 1(a) for an illustration).

\begin{figure}[ptb]
{\large \centering   \subfigure[WOS within the domain]{
\label{fig:subfig:a}     \includegraphics[width=0.37\textwidth]{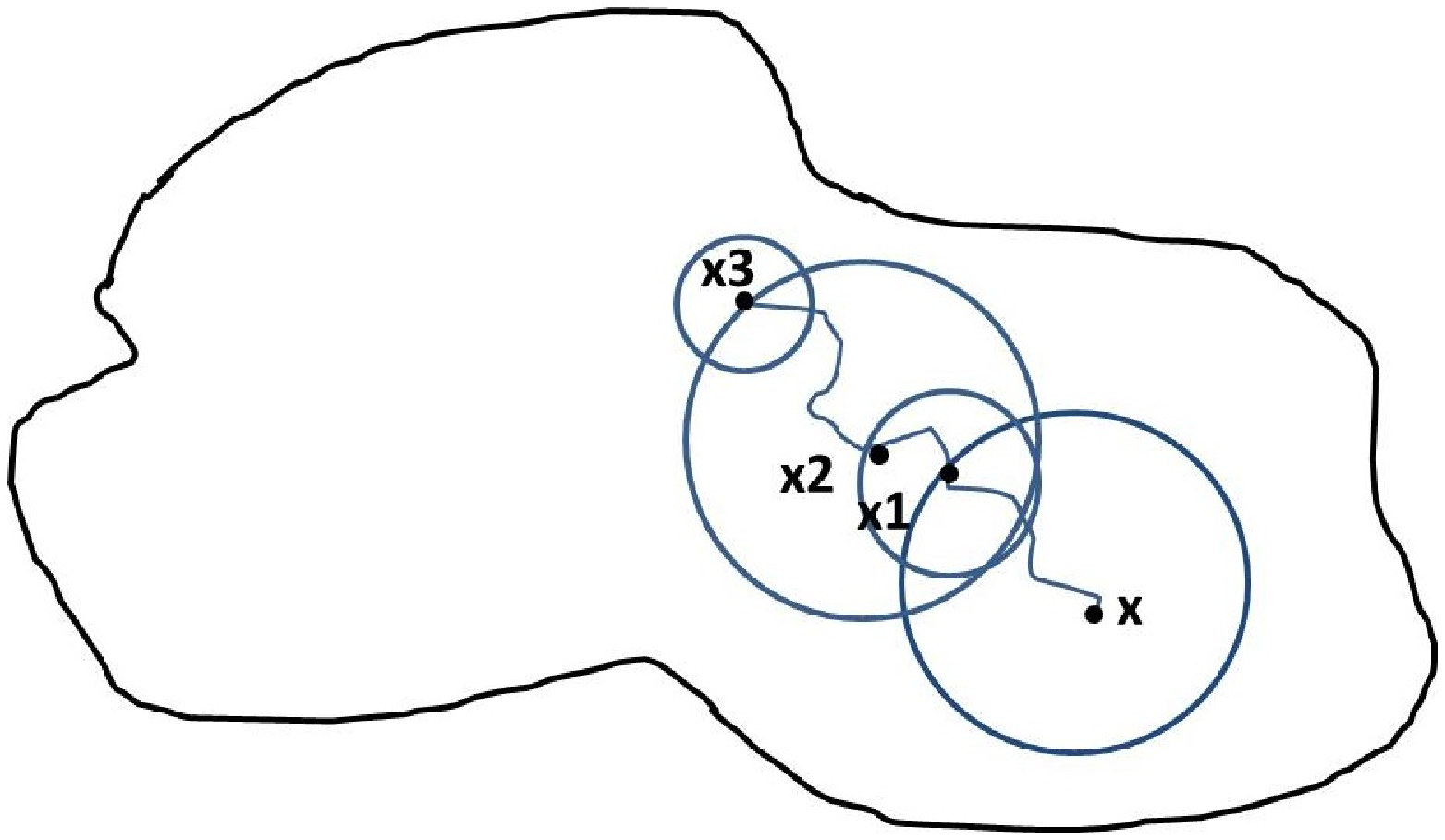}}
\hspace{1in}
\subfigure[WOS (with a maximal step size for each jump) within the domain]{
\label{fig:subfig:b}     \includegraphics[width=0.35\textwidth]{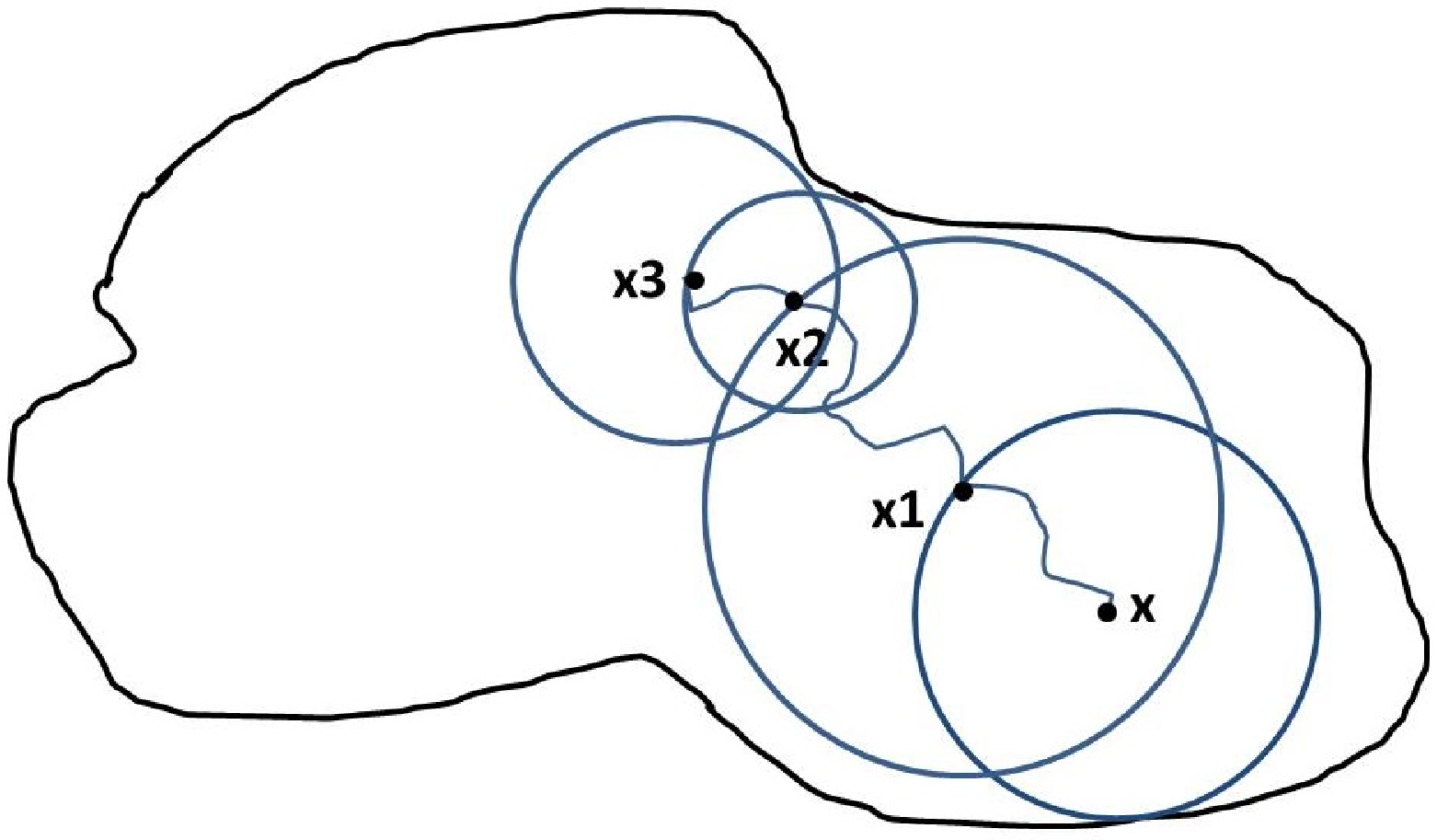}}
}\caption{Walk on Spheres method}%
\label{fig:subfig}%
\end{figure}

Now we can define an estimator of (\ref{eq13}) with $c=0$ by
\begin{equation}
u(x)\approx\frac{1}{N}\sum_{i=1}^{N}u(x_{i}),
\end{equation}
where $N$ is the number of Brownian paths sampled and $x_{i}$ is the first
hitting point of each path on the boundary. To speed up the WOS process,
maximum possible size of the sphere for each step would allow faster first
hitting on the boundary.

\begin{itemize}
\item \bigskip WOS and RBM
\end{itemize}

For the reflecting boundary, we will construct a strip region around the
boundary (see Fig. 2) and allow the process $X_{t}$ to move according to the
law of BM continuously. Before the path enters the strip region, the radius of
WOS is chosen to be of a maximum possible size less than the distance to the
boundary. Once the particle is in the strip region, the radius of the
WOS\ sphere is fixed at a constant $\Delta x$ (or $2\Delta x$, see Fig. 3).
With this approach, according to the definition ($\ref{eq5}$), the local time
may be interpreted as
\begin{equation}
dL(t)\approx\frac{\int_{t_{j-1}}^{t_{j}}I_{D_{\epsilon}}(X_{s})ds}{\epsilon},
\label{eq41}%
\end{equation}
which is
\begin{equation}
dL(t)\approx\frac{\int_{t_{j-1}}^{t_{j}}I_{D_{\epsilon}}(X_{s})ds}{\epsilon
}=(n_{t_{j}}-n_{t_{j-1}})\frac{(\Delta x)^{2}}{3\epsilon}, \label{eq43}%
\end{equation}
given a prefixed constant $\Delta x$ in the strip region and $n_{t_{j}}$ be
the cumulative steps that path stays within the $\epsilon$-region from the
begining until time $t_{j}$ (see Remark below for definition). Notice that
only those steps where the path of $X_{t}$ remains in the $\epsilon$-region
will contribute to $n_{t_{j}}$ because the SRBM may lie out of the $\epsilon
$-region at other steps. More details can be found in \cite{[4]}, where the
same construction is applied for the Neumann boundary value problem. One may
refer to Fig. 3 for an illustration of the behavior of path near the boundary.

\begin{figure}[ptb]
\centering {\large \includegraphics[width=0.38\textwidth]{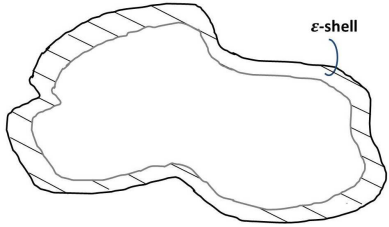}
}\caption{A $\epsilon$-region for a bounded domain in $R^{3}$}%
\end{figure}

\begin{figure}[ptb]
\centering {\large \includegraphics[width=0.5\textwidth]{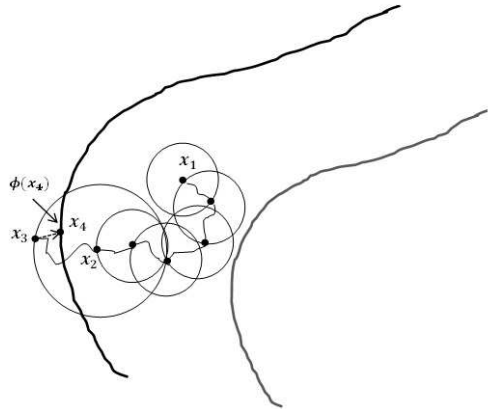}
}\caption{WOS in the $\epsilon$-region. At point $x_{1}$, BM path first hits
the $\epsilon$-region. By WOS with a prefixed radius $\Delta x$, the path
continues moving subsequently to $x_{2}$ where the distance to the boundary is
less than $\Delta x$. Enlarge the radius to $2\Delta x$, the path then have a
probability to run out of the domain to $x_{3}$. Pull back to the closest
point $x_{4}$ on the boundary, record $\phi(x_{4})$ and continue WOS-sampling
starting at $x_{4}$.}%
\end{figure}

\begin{rem}
Occupation time of SRBM $X_{t}$ in the numerator of ($\ref{eq41}$) was
calculated in terms of that of BM sampled by the walks on spheres. Notice here
that within the $\epsilon$-region, the radius of the WOS may be $\Delta x$ or
$2\Delta x$, which implies that the corresponding elapsed time of one step for
local time could be $(\Delta x)^{2}/3$ or $(2\Delta x)^{2}/3$. The latter is
four times bigger than the former. But if we absorb the factor $4$ into
$n_{t}$, $(\ref{eq43})$ still holds. In practical implementation, we treat
$n_{t}$ as a vector of entries of increasing value, the increment of each
component of $n_{t}$ over the previous one after each step of WOS will be 0, 1
or 4, corresponding to the scenarios that $X_{t}$ is out of the $\epsilon
$-region, in the $\epsilon$-region while sampled on the sphere of a radius
$\Delta x$, or in the $\epsilon$-region while sampled on the sphere of a
radius $2\Delta x$, respectively.
\end{rem}

Robin boundaries represent a general form of an insulating boundary condition
for convection-diffusion equations where $c(x)$ stands for the positive
diffusive coefficients. For our numerical test, we will consider two cases: a
positive constant $c$ and a positive function $c(x)$.

\subsection{Numerical Tests}

The numerical approximations obtained are compared to the true solutions on a
selected circle and a line segment, respectively, for the following three test
domains in $R^{3}$:

\begin{enumerate}
\item A cube centered at the origin with a length 2;

\item A sphere centered at the origin with a radius 1;

\item An ellipsoid centered at the origin with axial lengths [3, 2, 1].
\end{enumerate}

\bigskip

The location of the circle is given by%
\begin{equation}
\{(x,y,z)^{T}=(r\cos\theta_{1}\sin\theta_{2},r\sin\theta_{1}\sin\theta
_{2},r\cos\theta_{2})^{T}\}\label{eq45}%
\end{equation}
with $r=0.6$, $\theta_{1}=0:k\cdot2\pi/30:2\pi$, $\theta_{2}=\pi/4$ with
$k=1,...,15$. While the line segment is defined with endpoints
$(0.4,0.4,0.6)^{T}$ and $(0.1,0,0)^{T}$. Fifteen uniformly spaced points on
the line are selected to monitor the accuracy of the numerical solutions.

Finally, we set the true solution of the Robin boundary problem (\ref{eq17})
to be
\begin{equation}
u(x)=\sin3x\sin4y\ e^{5z}+5.\label{eq47}%
\end{equation}

\subsubsection{Constant $c(x)$}

\textbf{Example 1} \quad$c(X_{t})=1$

In this case, ($\ref{eq37}$) is reduced to
\begin{equation}
u(x)=E^{x}\{\int_{0}^{\infty}e^{\int_{0}^{t}dL_{t}}f(X_{t})dL_{t}\},
\label{eq49}%
\end{equation}
which is equivalent to
\begin{equation}
u(x)=E^{x}\{\int_{0}^{\infty}e^{L_{t}-L_{0}}f(X_{t})dL_{t}\} \label{eq51}%
\end{equation}
or
\begin{equation}
u(x)=E^{x}\{\int_{0}^{\infty}e^{L_{t}}f(X_{t})dL_{t}\}, \label{eq53}%
\end{equation}
for a starting point $x$ belonging to the interior of the solution domain.

We will truncate the time interval to $[0,T]$, an approximation to
($\ref{eq53}$) will be
\begin{equation}
\tilde{u}(x)=E^{x}\{\int_{0}^{T}e^{L_{t}}f(X_{t})dL_{t}\}.\label{eq55}%
\end{equation}
Using the fact that
\begin{equation}
dL_{t}\approx(n_{t}-n_{t-1})\frac{(\Delta x)^{2}}{3\epsilon},\label{eq57}%
\end{equation}

we can rewrite ($\ref{eq55}$) as
\begin{equation}
\tilde{u}(x)=E^{x}\{\int_{0}^{T}e^{n_{t}\frac{(\Delta x)^{2}}{3\epsilon}%
}f(X_{t})(n_{t}-n_{t-1})\frac{(\Delta x)^{2}}{3\epsilon}\}.
\end{equation}

Next identifying the time interval with the length of sample path NP, we have
\begin{equation}
\tilde{u}(x)=E^{x}\left\{  \sum_{j^{\prime}=0}^{NP}e^{n_{t_{j}}\frac{(\Delta
x)^{2}}{3\epsilon}}f(X_{t_{j}})(n_{t_{j}}-n_{t_{j-1}})\frac{(\Delta x)^{2}%
}{3\epsilon}\right\}  ,\label{eq59}%
\end{equation}
where $j^{\prime}$ denotes each step of the path and $j$ denotes the steps
where the path hits the boundary.

At each step along a path we first evaluate
\[
e^{n_{t_{j}}\frac{(\Delta x)^{2}}{3\epsilon}}f(X_{t_{j}})(n_{t_{j}}%
-n_{t_{j-1}})\frac{(\Delta x)^{2}}{3\epsilon},
\]
if $X_{t_{j}}$ hits the boundary, we then compute $f(X_{t_{j}})(n_{t_{j}%
}-n_{t_{j-1}})\frac{(\Delta x)^{2}}{3\epsilon}$, followed by multiplying it by
$e^{n_{t_{j}}\frac{(\Delta x)^{2}}{3\epsilon}}$, which uses the cumulative
time of $L_{t_{j}}$ from $t=0$ to $t_{j}$. Finally, the expectation is done
via the average over $N$ sample paths.

The simulation results of a cubic domain are presented in Fig. 4 and 5. The
two figures show the convergency of the approximations as the length of path
increases from $1.35e4$ to $1.43e4$ and $1.6e4$ to $1.7e4$ over the circle and
the line segment, respectively. Some deviations are seen at the tail in Figure
4(a) and among the middle points in Figure 4(b). Meanwwhile, for the spherical
and ellipsoid domains (Figure 6 and 7), the approximations are better and the
errors are relatively smaller especially over the line segments, which are
below 3\% in Figure 6(b) and Figure 7(b).

\begin{figure}[ptb]
{\large \centering   \subfigure[$\epsilon$ = 3$\Delta x$, Err = 9.59\%, NP=1.35e4, $\Delta x$=5e-4]{
\label{fig:subfig:a}
\includegraphics[width=0.38\textwidth]{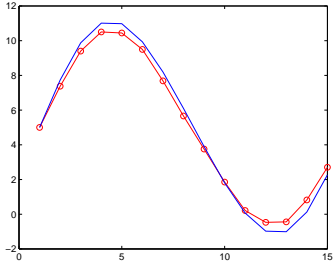}} \hspace{1in}
\subfigure[$\epsilon$ = 4$\Delta x$, Err =8.84\%, NP=1.6e4, $\Delta x$=5e-4]{
\label{fig:subfig:b}
\includegraphics[width=0.38\textwidth]{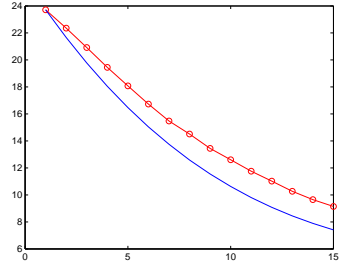}} }\caption{Cubic domain:
number of paths $N=2e5$ and $c(X_{t})=1$. (Left - circle; right - line
segement) }%
\label{fig:subfig}%
\end{figure}

\begin{figure}[ptb]
{\large \centering   \subfigure[$\epsilon$ = 3$\Delta x$, Err = 5.50\%, NP=1.43e4, $\Delta x$=5e-4]{
\label{fig:subfig:a}
\includegraphics[width=0.38\textwidth]{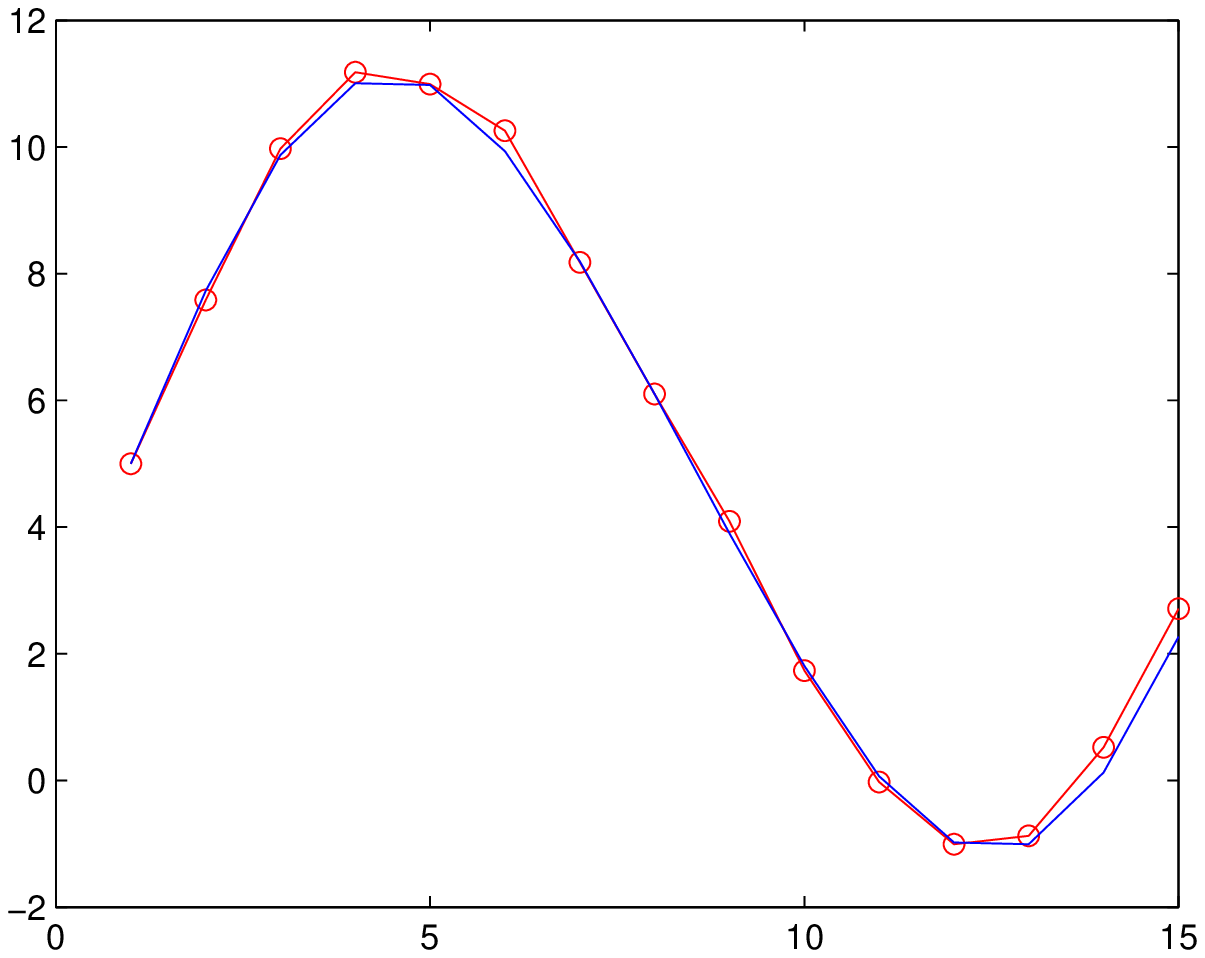}} \hspace{1in}
\subfigure[$\epsilon$ = 4$\Delta x$, Err = 6.49\%, NP=1.7e4, $\Delta x$=5e-4]{
\label{fig:subfig:b}
\includegraphics[width=0.38\textwidth]{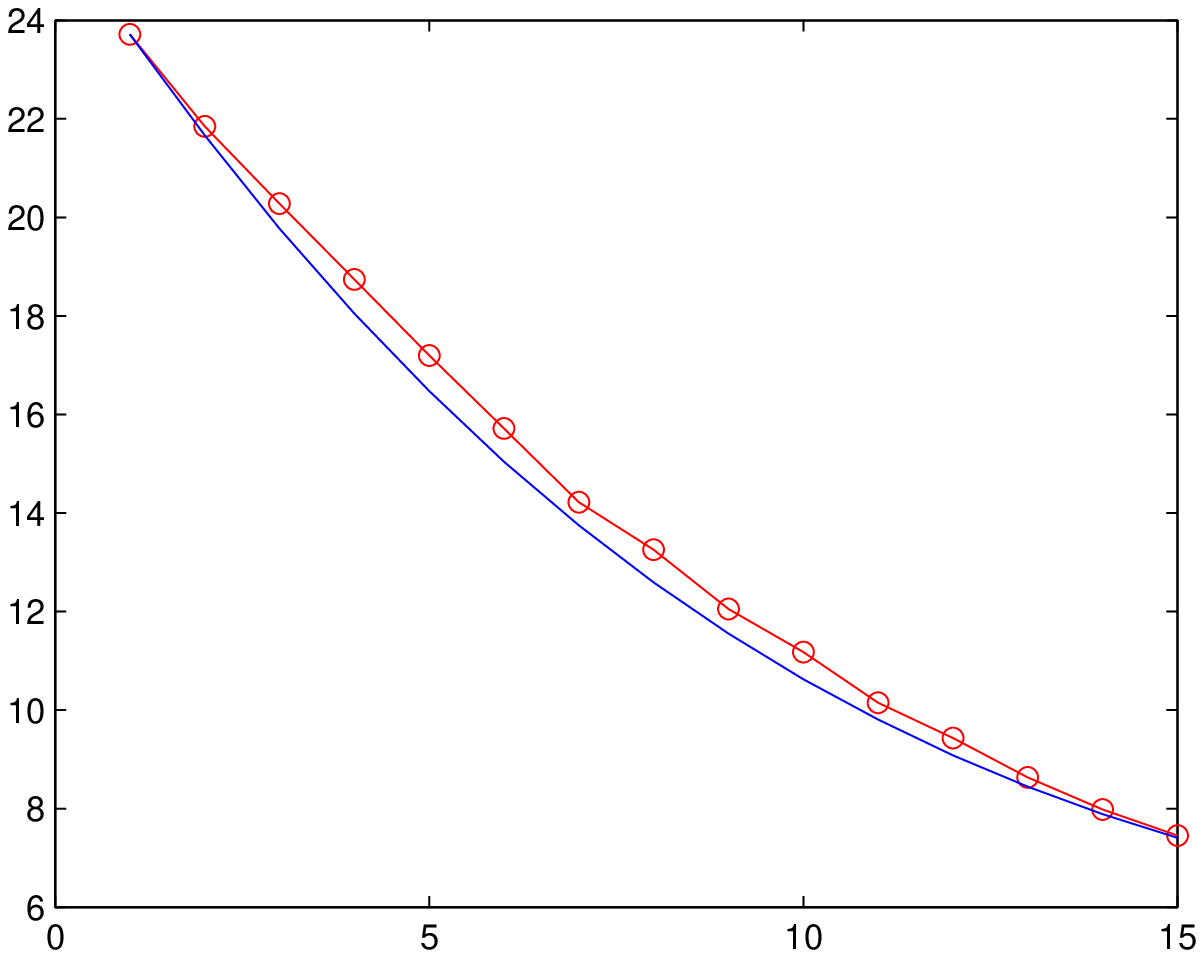}} }\caption{Cubic domain:
number of paths $N=2e5$ and $c(X_{t})=1$. (Left - circle; right - line
segement) }%
\label{fig:subfig}%
\end{figure}

\begin{figure}[ptb]
{\large \centering   \subfigure[$\epsilon$ = 3$\Delta x$, Err = 3.96\%, NP=6e3, $\Delta x$=5e-4]{
\label{fig:subfig:a}
\includegraphics[width=0.38\textwidth]{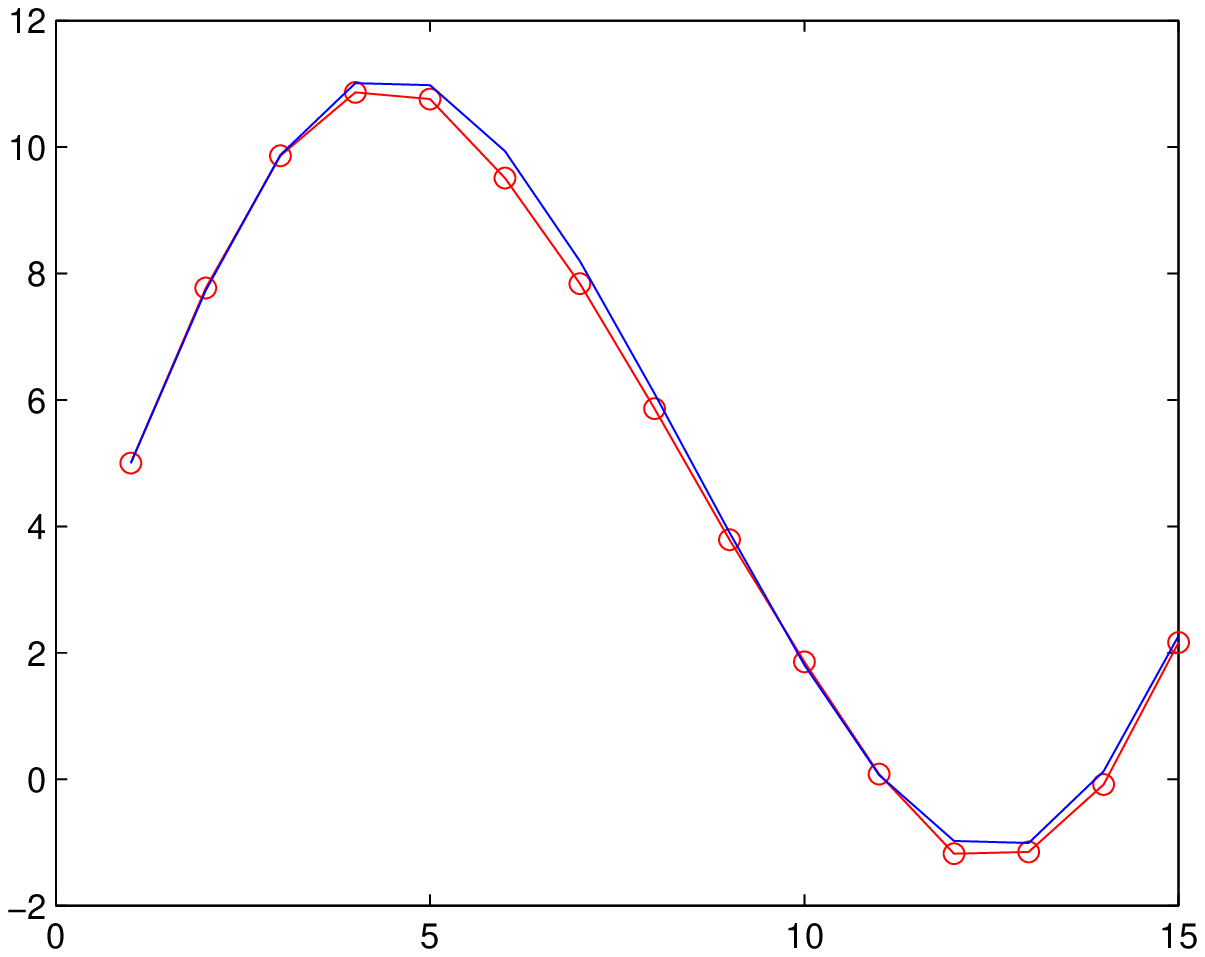}} \hspace{1in}
\subfigure[$\epsilon$ = 3$\Delta x$, Err = 1.24\%, NP=5.5e3, $\Delta x$=5e-4]{
\label{fig:subfig:b}
\includegraphics[width=0.38\textwidth]{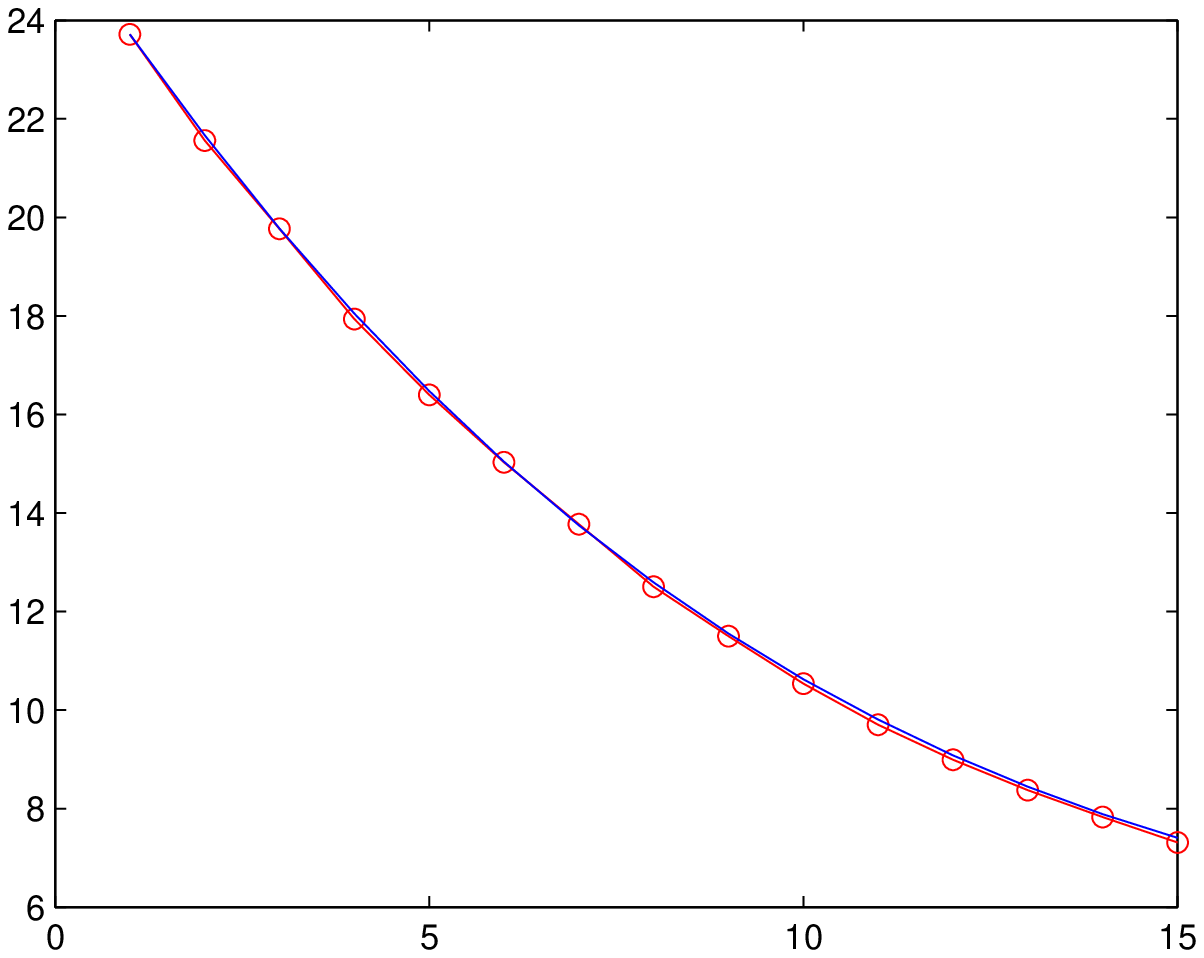}} }\caption{Spherical domain:
number of paths $N=2e5$ and $c(X_{t})=1$. (Left - circle; right - line
segement)}%
\label{fig:subfig}%
\end{figure}

\begin{figure}[ptb]
{\large \centering   \subfigure[$\epsilon$ = 3$\Delta x$, Err = 2.42 \%, NP= 4.5e3, $\Delta x$=5e-4]{
\label{fig:subfig:a}
\includegraphics[width=0.38\textwidth]{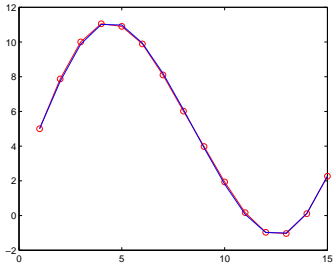}} \hspace{1in}
\subfigure[$\epsilon$ = 3$\Delta x$, Err = 2.44\%, NP=4.5e3, $\Delta x$=5e-4]{
\label{fig:subfig:b}
\includegraphics[width=0.38\textwidth]{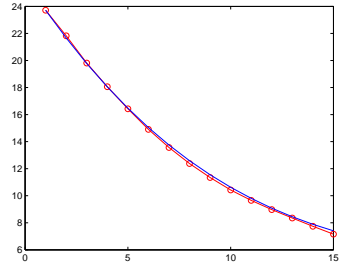}} }\caption{Ellipsoid domain:
number of paths $N=2e5$ and $c(X_{t})=1$. (Left - circle; right - line
segement)}%
\label{fig:subfig}%
\end{figure}

\subsubsection{Variable c(x)}

\textbf{Example 2}\ $c(X_{t})=|x|$, $x$ is the first component of $X_{t}$ on
the boundary. Similar to Example 1, we have%

\begin{equation}
u(x)=E^{x}\left\{  \int_{0}^{\infty}e^{\int_{0}^{t}c(X_{s})dL_{s}}%
f(X_{t})dL_{t}\right\}  .\label{eq61}%
\end{equation}
It can be seen that $c(X_{s})dL_{s}$ and $f(X_{t})dL_{t}$ have the same form,
so we can handle $c(X_{s})dL_{s}$ exactly the same way as $f(X_{t})dL_{t}$.
Then, we have%

\begin{equation}
u(x)=E^{x}\left\{  \sum_{j^{\prime}=0}^{NP}e^{\sum_{k=0}^{j}c(X_{t_{k}%
})(n_{t_{k}}-n_{t_{k-1}})\frac{h^{2}}{3\epsilon}}f(X_{t_{j}})(n_{t_{j}%
}-n_{t_{j-1}})\frac{h^{2}}{3\epsilon}\right\}  . \label{eq63}%
\end{equation}

Notice that the term
\begin{equation}
e^{\sum_{k=0}^{j}c(X_{t_{k}})(n_{t_{k}}-n_{t_{k-1}})\frac{h^{2}}{3\epsilon}%
}\label{eq65}%
\end{equation}
cumulates all the information of $c(X_{t})$ with respect to the local time
from the beginning to the current time. If $c(X_{t})=|x|$, then%

\begin{equation}
u(x)=E^{x}\left\{  \sum_{j^{\prime}=0}^{NP}e^{\sum_{k=0}^{j}|x_{t_{k}}|%
(n_{t_{k}}-n_{t_{k-1}})\frac{h^{2}}{3\epsilon}}f(X_{t_{j}})(n_{t_{j}%
}-n_{t_{j-1}})\frac{h^{2}}{3\epsilon}\right\}  , \label{eq67}%
\end{equation}
where $j^{\prime}$ denote each step for the path and $j$ denotes the steps
where the path hits the boundary.

Numerical results are shown in Figure 8-10 for a cubic, a spherical and an
ellipsoid domain, respectively with some adjustment in $\Delta x$ and $NP$.
Here we still have similar results for cube with errors around 6.5\%. For the
sphere, we change $\Delta x$ to $4e-4$ and there are deviation around the
middle in Figure 9(a) which may explain the overall error only 6.74\% while it
performs well over the line segment in Figure 9(b) with a smaller error of
3.1\%. For the ellipsoid, the results are similar as in Example 1 and maintain
an error below 4\%.

\begin{figure}[ptb]
{\large \centering   \subfigure[$\epsilon$ = 3$\Delta x$, Err = 6.35\%, NP=1.6e4, $\Delta x$=5e-4]{
\label{fig:subfig:a}
\includegraphics[width=0.38\textwidth]{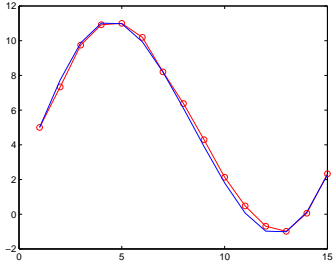}} \hspace{1in}
\subfigure[$\epsilon$ = 3$\Delta x$, Err = 6.66\%, NP=1.48e4, $\Delta x$=5e-4]{
\label{fig:subfig:b}
\includegraphics[width=0.38\textwidth]{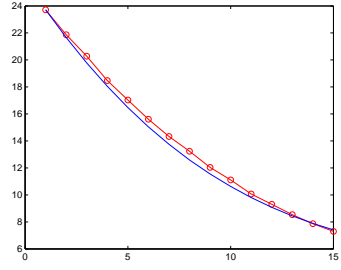}} }\caption{Cubic domain:
number of paths $N=2e5$ and $c(X_{t})=|x|$. (Left - circle; right - line
segement) }%
\label{fig:subfig}%
\end{figure}

\begin{figure}[ptb]
{\large \centering   \subfigure[$\epsilon$ = 3$\Delta x$, Err = 6.74\%, NP=6.5e3, $\Delta x$=4e-4]{
\label{fig:subfig:a}
\includegraphics[width=0.38\textwidth]{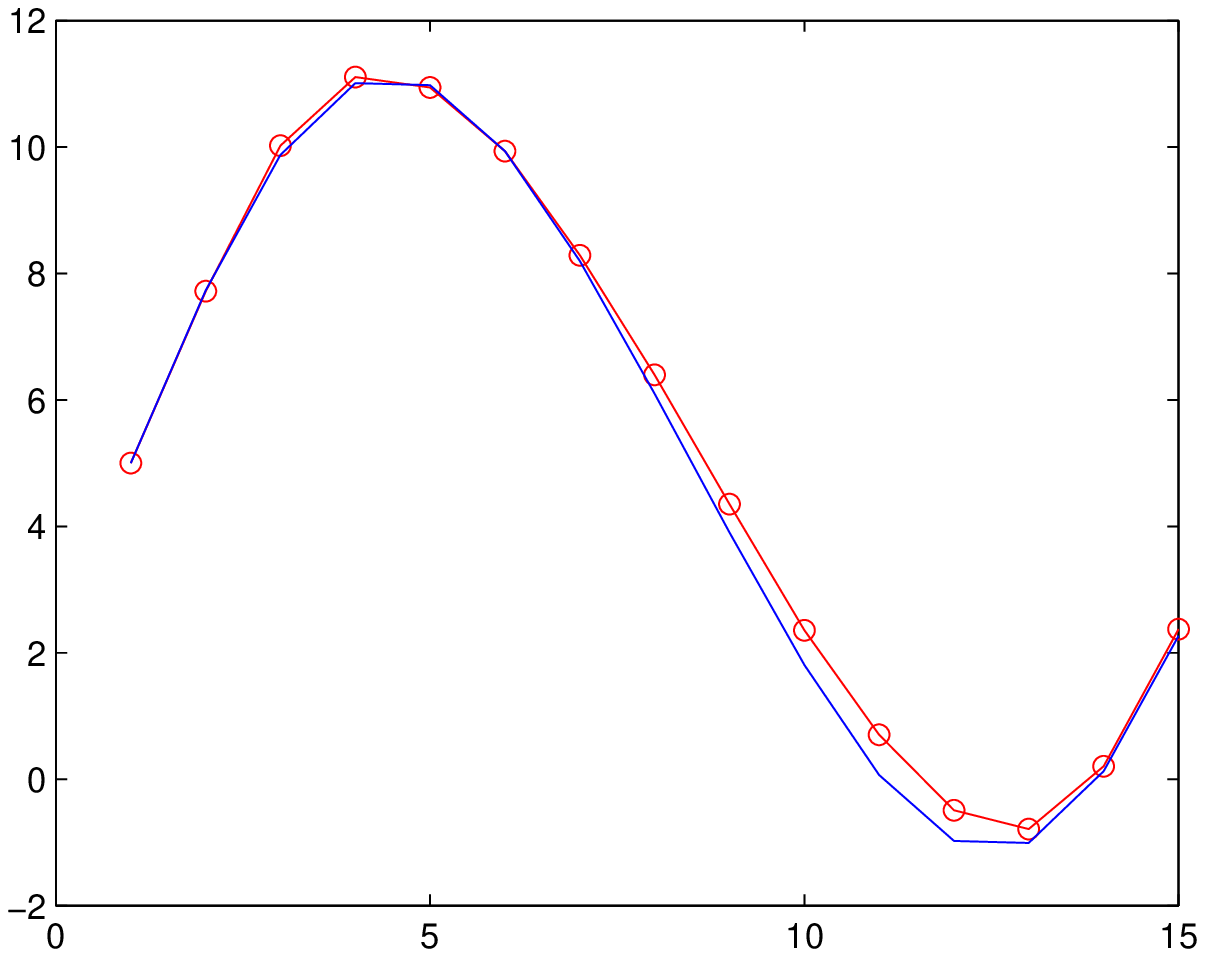}} \hspace{1in}
\subfigure[$\epsilon$ = 3$\Delta x$, Err = 3.10\%, NP=6e3, $\Delta x$=4e-4]{
\label{fig:subfig:b}
\includegraphics[width=0.38\textwidth]{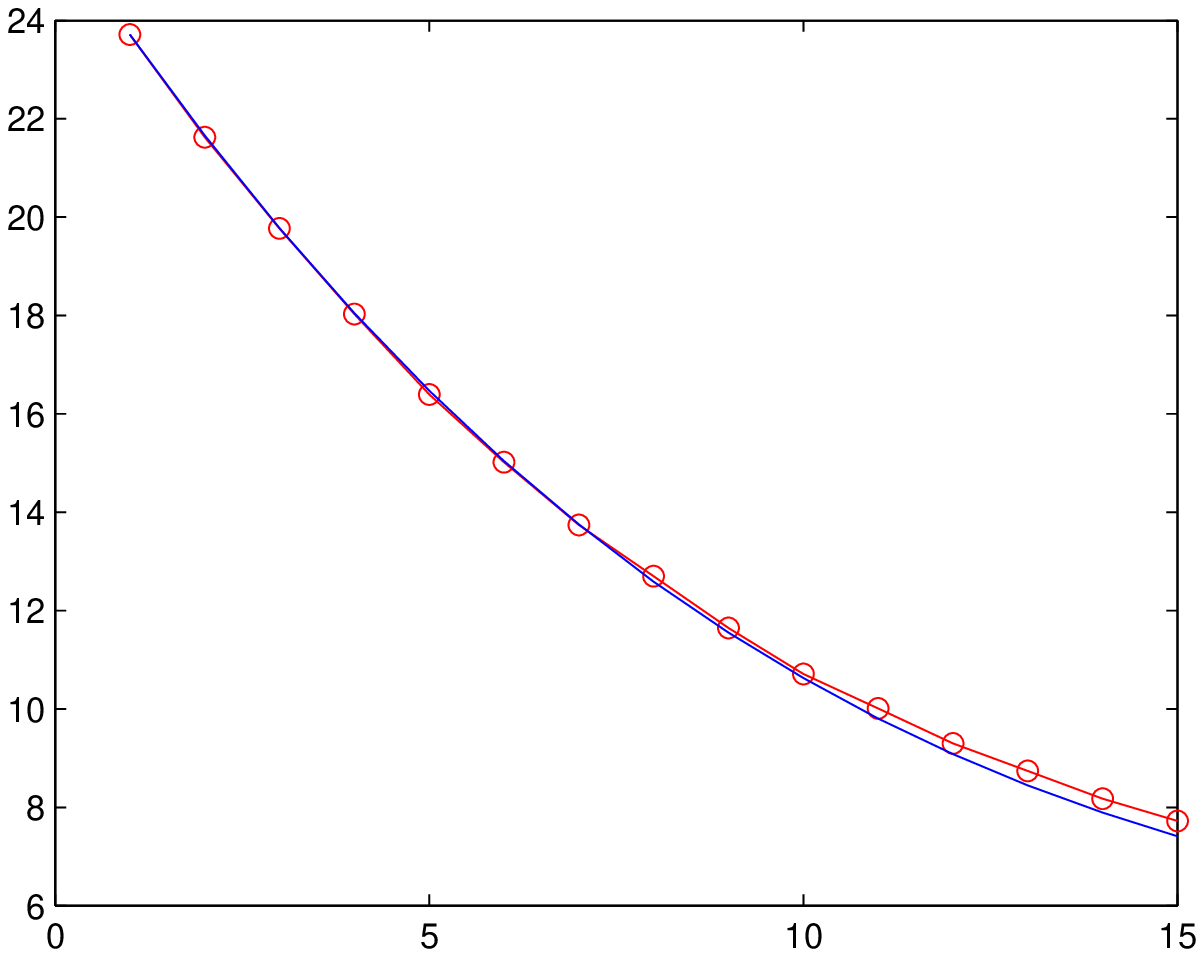}} }\caption{Spherical domain:
number of paths $N=2e5$ and $c(X_{t})=|x|$. (Left - circle; right - line
segement)}%
\label{fig:subfig}%
\end{figure}

\begin{figure}[ptb]
{\large \centering   \subfigure[$\epsilon$ = 3$\Delta x$, Err = 3.76 \%, NP= 5e3, $\Delta x$=4e-4]{
\label{fig:subfig:a}
\includegraphics[width=0.38\textwidth]{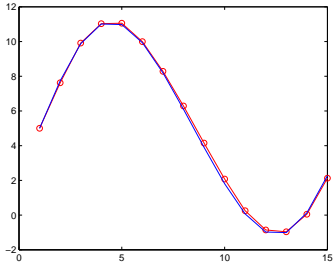}} \hspace{1in}
\subfigure[$\epsilon$ = 3$\Delta x$, Err = 1.93\%, NP=5e3, $\Delta x$=4e-4]{
\label{fig:subfig:b}
\includegraphics[width=0.38\textwidth]{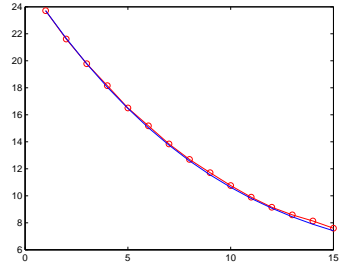}} }\caption{Ellipsoid domain:
number of paths $N=2e5$ and $c(X_{t})=|x|$. (Left - circle; right - line
segement)}%
\label{fig:subfig}%
\end{figure}

\section{Conclusions and future work}

This paper presents a Monte Carlo simulation method to solve the third
boundary problems associated with Laplace equations. The idea of simulating
sample paths of SRBM by the WOS within the strip region shows its efficiency
and accuracy in estimating local time and evaluating Feynman-Kac formula. It
should be noted that the cases that $q\neq0$ needs further work due to the
unknown exit time out of the sphere at each step. For the Poisson equation,
the contribution of the source term might be computed as a conditional
integral \cite{[19]}. Moreover, the proper truncation of time period is
unknown, though it is proven that the variance of the approximation increases
linearly of $T$ \cite{[13]}.

For future work, more flexible domains with local convexity will be considered
as it relates to the calculations of electrical properties such as the
conductivity of composite materials where the particle shapes plays an
important role \cite{[20]}.

\section*{Acknowledgement}

The authors Y.J.Z and W.C. acknowledge the support of the National Science
Foundation (DMS-1315128) and the National Natural Science Foundation of China
(No. 91330110) for the work in this paper.

\end{document}